\documentclass[11pt]{amsart}
\usepackage{amsmath,amsthm,amsfonts,amscd,amssymb,eucal,latexsym}
\usepackage[all]{xy}

\theoremstyle{definition}
\newtheorem{theorem}{Theorem}[section]
\newtheorem{corollary}[theorem]{Corollary}
\newtheorem{lemma}[theorem]{Lemma}
\newtheorem{proposition}[theorem]{Proposition}

\newtheorem{question}[theorem]{Question}
\newtheorem{problem}[theorem]{Problem}

\newcommand{\topol}{{\rm top}}
\newcommand{\rank}{{\rm rank}}
\newcommand{\CPA}{{\rm CPA}}
\newcommand{\rcp}{{\rm rcp}}

\newcommand{\sep}{{\rm sep}}
\newcommand{\spa}{{\rm spn}}
\newcommand{\spn}{{\rm span}}

\newcommand{\id}{{\rm id}}

\allowdisplaybreaks

\begin{document}

\title[A geometric approach to Voiculescu-Brown entropy]{A geometric 
approach to Voiculescu-Brown entropy}

\author{David Kerr}
\address{Dipartimento di Matematica, Universit\`{a} di Roma 
``La Sapienza,'' P.le Aldo Moro, 2, 00185 Rome, Italy}
\email{kerr@mat.uniroma1.it}
\date{April 28, 2003}

\begin{abstract}
A basic problem in dynamics is to identify systems 
with positive entropy, i.e., systems which are ``chaotic.'' While
there is a vast collection of results addressing this issue in 
topological dynamics, the phenomenon of positive entropy remains by and 
large a mystery within the broader noncommutative domain of $C^*$-algebraic 
dynamics. To shed some light on the noncommutative situation we propose 
a geometric perspective inspired by work of Glasner and Weiss on 
topological entropy. 
This is a written version of the author's talk
at the Winter 2002 Meeting of the Canadian Mathematical Society
in Ottawa, Ontario.
\end{abstract}

\maketitle

\section{Introduction}

Entropy has been very successful in 
topological dynamics (as in measurable dynamics) as a numerical invariant  
measuring the complexity of a dynamical system (see, e.g., \cite{DGS,HK}). 
Topological entropy was introduced by Adler, Konheim, and McAndrew in 
\cite{AKM} with a definition based on open covers. Equivalent definitions 
in terms of separated and spanning sets with respect to a metric were given 
by Bowen \cite{Bow} and Dinaburg \cite{Din}. Much more recently Voiculescu
introduced a notion of entropy for automorphisms of unital nuclear 
$C^*$-algebras based on local approximation \cite{Voi}, and Brown
subsequently extended this to automorphisms of exact $C^*$-algebras
using nuclear embeddability \cite{Br}. By 
\cite[Prop.\ 4.8]{Voi} the topological entropy of a homeomorphism of a
compact metric space coincides with the Voiculescu-Brown entropy of the 
induced automorphism of the $C^*$-algebra of functions on the space,
and so Voiculescu-Brown entropy is an extension of 
topological entropy to the noncommutative domain (indeed Voiculescu and
Brown refer to their invariants as ``topological entropy,'' but we have 
refrained from this terminology to avoid confusion). However, because 
Voiculescu's idea of using local approximation constitutes a
fundamentally reconceptualized approach to defining entropy (as
necessitated by its analytic context,
which does not allow for any kind of direct analogue of an open cover with
which a genuine dynamical invariant may be obtained), Voiculescu-Brown
entropy ultimately requires tools of a completely different 
formal and technical nature for its study.

To date this study has focused on
three aspects: computations for canonical examples (see, e.g., 
\cite{ECCE,BG,TFP,DS}), 
behaviour under taking crossed products \cite{Br,PS} and reduced free products 
\cite{BDS}, and the variational principle in the presence of a strategic 
amount of commutativity \cite{VP,PE,KP}. See \cite{NE} for a
survey. The methods for obtaining non-zero lower bounds in
computations have been invariably rooted in classical considerations. Indeed
they have involved either (i) relating the system to a
topological dynamical system which is a priori known to have positive
topological entropy, or (ii) appealing to measure-theoretic invariants like
Connes-Narnhofer-Thirring or Sauvageot-Thouvenot entropy, which themselves
use classical measurable partitions in their definitions. Thus the
fundamental problem of identifying systems with positive entropy
(i.e., systems which exhibit ``chaotic'' behaviour),
in however specific or general a setting, has not been addressed from
a broad noncommutative viewpoint, and so the phenomenon of positive entropy
has remained a mystery beyond the scope of commutativity (and even there 
a considerable degree of mystery persists, as we do not have a 
``topological'' proof of the equality of Voiculescu-Brown entropy and the 
topological entropy of the induced homeomorphism on the pure state space
in the separable unital commutative case,
which was established in \cite[Prop.\ 4.8]{Voi} using the
classical variational principle and continuity properties of 
Connes-Narnhofer-Thirring entropy).
As a consequence Voiculescu-Brown entropy has played a rather isolated 
role in $C^*$-dynamics, in contrast to the pervasive 
presence of topological entropy in topological dynamics. 

In fact the theory 
of $C^*$-dynamical systems has been much more concerned with questions of
$C^*$-algebraic structure and classification than with the strict 
investigation of long-term behaviour that is at the heart of topological
dynamics and for which entropy is a key tool. Symptomatic of the difficulties
in analyzing the long-term behaviour of noncommutative systems is the
typical lack of discrete data that can be workably assembled to yield
a meaningful statement about some aspect of the dynamics.
This is starkly illustrated by the collection of automorphisms of the rotation 
$C^*$-algebras $A_\theta$ associated to a given matrix in $SL(2,\mathbb{Z})$
with eigenvalues off the unit circle: for a subset of $\theta\in [0,1)$ of 
full Lebesgue measure the canonical tracial state is the 
unique invariant state \cite{NT}, while in 
the ``degenerate'' case $\theta = 0$, in which we recover the corresponding
hyperbolic toral automorphism on the pure state space, there is a rich 
supply of periodic points in terms of which much about the system can
be expressed, including the entropic growth (see, e.g., \cite{HK}). We
can also make a comparison here with the viewpoint of semiclassical 
analysis, which, as opposed to directly extending topological dynamics
to incorporate the noncommutative case, 
extracts discrete spectral data from within 
the matrix framework associated to the values $\theta = 1/n$ 
(``quantization'') and performs an asymptotic analysis thereupon,
with the possibility of establishing a correspondence with classical
information in the limit $n\to\infty$. See the introduction to \cite{Ze}
for a discussion and \cite{KR1,KR2,KR3} for some recent results.

The purpose of this article is to introduce a geometric perspective
that yields some insight into the mechanisms behind positive 
Voiculescu-Brown entropy. Our inspiration lies in the link between 
topological dynamics and the geometric theory of Banach spaces that was 
established by Glasner and Weiss in one of the two proofs they gave in 
\cite{GW} for the following striking result.

\begin{theorem}[{\cite[Thm.\ A]{GW}}]\label{T-GW}
If a homeomorphism from a compact metric space $X$ to itself has zero 
topological entropy, then so does the induced homeomorphism on the space 
of probability measures on $X$.
\end{theorem}

\noindent Since the Voiculescu-Brown entropy agrees with the topological 
entropy on the pure state space for automorphisms of separable unital 
commutative $C^*$-algebras, Glasner and Weiss's 
result is equivalent to the assertion that zero Voiculescu-Brown entropy 
for such an automorphism implies zero topological entropy for the induced
homeomorphism on the state 
space. We have shown that this assertion in fact holds for automorphisms
of any separable unital exact $C^*$-algebra, and we can additionally drop 
the assumption of a unit by replacing the state space with the quasi-state 
space in the general exact setting (Theorem~\ref{T-zero}). 
Full details of this result can be found in \cite{EID}. The key 
geometric tool is the asymptotic exponential 
dependence of $k$ on $n$ given an approximately isometric embedding
of $\ell^n_1$ into the $C^*$-algebra $M_k$ of $k\times k$ matrices. This 
geometric fact also has the consequence that the presence of a 
certain supply of dynamically generated subspaces approximately isometric 
to $\ell^n_1$ is sufficient to obtain positive Voiculescu-Brown entropy
(Proposition~\ref{P-pdpe}). In particular 
we can show, without relying in any way on classical dynamical 
concepts, that certain $C^*$-dynamical systems constructed in an 
operator-theoretic manner have positive entropy. We also obtain some
information about the behaviour of Voiculescu-Brown entropy under taking
extensions.

The main body of the paper consists of three sections.
Sections \ref{S-ell1} and \ref{S-EID} revolve around the results described
above involving positive entropy and
approximately isometric embeddings of $\ell^n_1$ into $M_k$
at the $C^*$-algebra and state space levels, respectively, while in 
Section~\ref{S-free} we apply our geometric perspective in a complementary  
way with a look at the free shift on 
$C^*_r (\mathbb{F}_\infty )$ as an example of subexponential dynamical growth.
\medskip

\noindent{\it Acknowledgements.} 
This work was supported by the Natural 
Sciences and Engineering Research Council of Canada and was carried out 
during stays at the University of Tokyo and the University of Rome 
``La Sapienza'' over the academic years 2001--2002 and 2002--2003,
respectively. 
I thank Yasuyuki Kawahigashi at the University of Tokyo
and Claudia Pinzari at the University of Rome ``La Sapienza'' for their
generous hospitality. 
I would also like to thank the Canadian Mathematical
Society for the opportunity to present this material at the
Winter 2002 Meeting.

\section{Entropy and embeddings of $\ell^n_1$ into $M_k$}\label{S-ell1}

We begin by recalling the definitions of topological entropy and
Voiculescu-Brown entropy. Let $X$ be a compact metric space and $T : X\to X$
a homeomorphism. For a finite open cover $\mathcal{U}$ we set
$$ h_{\topol}(T, \mathcal{U}) = \lim_{n\to\infty} \frac1n \log
N(\mathcal{U}\vee T^{-1}\mathcal{U}\vee\cdots\vee T^{-(n-1)} \mathcal{U}) $$
where $N(\cdot )$ denotes the smallest cardinality of a subcover and
the join $\mathcal{U}_1 \vee\cdots\vee \mathcal{U}_m$ of a finite collection
$\mathcal{U}_1 , \dots , \mathcal{U}_m$ of open covers is the
set of all non-empty intersections $U_1 \cap\cdots\cap U_m$ with
$U_i \in\mathcal{U}_i$ for each $i=1 , \dots , m$.
The {\it topological entropy} of $T$ is defined by
$$ h_{\topol}(T) = \sup_{\mathcal{U}} h_{\topol}(T, \mathcal{U}) $$
where the supremum is taken over all open covers $\mathcal{U}$. We can
alternatively express the entropy in terms of separated and spanning sets.
A set $E\subset X$ is said to be 
{\it $(n,\varepsilon )$-separated (with respect to $T$)} if for every $x,y\in 
E$ with $x\neq y$ there is a $0\leq k \leq n-1$ such that $d(T^k x , T^k y)
> \varepsilon$, and {\it $(n,\varepsilon )$-spanning (with respect to $T$)} 
if for every $x\in X$ there is a $y\in E$ such that $d(T^k x , T^k y ) \leq 
\varepsilon$ for each $k=0, \dots , n-1$. We write $\sep_n (T,\varepsilon )$
and $\spa_n (T,\varepsilon )$ to denote the largest 
cardinality of an $(n,\varepsilon )$-separated set and the smallest 
cardinality of an $(n,\varepsilon )$-spanning set, respectively. We then have
$$ h_{\topol}(T) = \sup_{\varepsilon >0} \limsup_{n\to\infty}\frac1n
\log\sep_n (T,\varepsilon ) = \sup_{\varepsilon >0} \limsup_{n\to\infty}\frac1n
\log\spa_n (T,\varepsilon ) . $$ The fundamental prototypical example is
the shift on $\{ 1, \dots , d \}^{\mathbb{Z}}$, with entropy $\log d$.
For general references on topological entropy see \cite{DGS,HK}.

Turning now to the noncommutative domain, we let $A$ be an exact 
$C^*$-algebra and $\pi : A \to\mathcal{B}(\mathcal{H})$
a faithful representation. By \cite{Kir} exactness is equivalent to nuclear
embeddability, and the latter guarantees, for every
finite $\Omega\subset A$ and $\delta > 0$, the non-emptiness of the 
collection $\CPA (\pi , \Omega , \delta )$ 
of triples $(\phi , \psi , B)$ where $B$ is a 
finite-dimensional $C^*$-algebra and $\phi : A\to B$ and $\psi : B\to
\mathcal{B}(\mathcal{H})$ are completely positive contractive linear maps 
such that the diagram
\begin{gather*}
\xymatrix{
A \ar[dr]_{\phi} 
\ar[rr]^{\pi} & & \mathcal{B} (\mathcal{H}) \\
& B \ar[ur]_{\psi}}
\end{gather*}
approximately commutes to within $\delta$ on $\Omega$, i.e., 
$\| (\psi\circ\phi )(x) - \pi (x) \| < \delta$ for all $x\in\Omega$.
As shown in the proof of \cite[Prop.\ 1.3]{Br}, the infimum of 
$\rank\, B$ over all $(\phi , \psi , B)\in\CPA (\pi , \Omega , \delta )$ 
(with rank referring to the dimension of a maximal commutative
$C^*$-subalgebra) does not depend on the particular faithful 
representation $\pi$; we denote this 
quantity by $\rcp (\Omega , \delta )$. We point out that the above 
$C^*$-algebras $B$ may in fact be taken to be full matrix algebras,
since a finite-dimensional $C^*$-algebra $B$ can be embedded in the
matrix algebra $M_k$ where $k = \rank\, B$, in which case the 
identity map on $B$ factors as the composition of the inclusion 
$B\hookrightarrow M_k$ with a conditional expectation $M_k \to B$. 
We also note that 
if $A$ is nuclear (in particular, if $A$ is commutative) we can alternatively 
define $\rcp (\Omega , \delta )$ by substituting the identity map on $A$ 
for $\pi$ in the above and taking the corresponding infimum. Furthermore,
if $A$ is unital then in both this nuclear reformulation and the more 
general exact setting we will obtain the same value of entropy
in the last line of the display below if we 
assume that the $\phi$ and $\psi$ in our 
approximately commuting diagrams are unital completely positive maps,
although $\rcp (\Omega , \delta )$ as we have defined it may not always be 
equal to its unital version
(see \cite{Br}). For an automorphism $\alpha$ of $A$ we then set 
\begin{align*}
ht(\alpha , \Omega , \delta ) &= \limsup_{n\to\infty}\frac1n \log\rcp 
(\Omega\cup\alpha\Omega\cup\cdots\cup\alpha^{n-1}\Omega , \delta ) , \\
ht(\alpha , \Omega ) &= \sup_{\delta > 0}ht(\alpha , \Omega , 
\delta ) , \\
ht(\alpha ) &= \sup_\Omega ht(\alpha , \Omega ) ,
\end{align*}
with the last supremum taken over all finite sets $\Omega\subset A$.
We call $ht(\alpha )$ the {\it Voiculescu-Brown entropy} of $\alpha$.
For some computations see \cite{Voi,ECCE,BG,TFP,DS} and for a survey see 
\cite{NE}.

Postponing momentarily the presentation of examples (to which we will turn
at natural points as our discussion evolves), we first make the general
remark that it is not at all clear from the definition
of Voiculescu-Brown entropy how exponential growth (i.e., positive values) 
might be produced, even for commutative systems. This is in
striking contrast to topological entropy, for which the mechanism behind
exponential growth is manifest in an example like the shift on
$\{ 1, \dots , d \}^{\mathbb{Z}}$.

In fact, in every case to date in which positive Voiculescu-Brown entropy
has been established, the argument has ultimately hinged on the use of some
measure-theoretic entropy, whether it has involved an appeal to 
Connes-Narnhofer-Thirring entropy or Sauvageot-Thouvenot entropy or to
the equality of Voiculescu-Brown entropy and the topological entropy on
the pure state space in the separable unital commutative setting (whose only
known proof relies on properties of Connes-Narnhofer-Thirring entropy and 
the classical variational principle \cite[Prop.\ 4.8]{Voi}). It is thus
natural to ask if we can avoid measure-theoretic entropies 
altogether and obtain a direct geometric explanation for the production
of exponential growth.

To approach this problem, let us first examine how the kind of mixing that
produces positive topological entropy is reflected geometrically at the
$C^*$-algebra level. Consider the right shift $T$ on 
$X = \{ -1,1 \}^{\mathbb{Z}}$ and the associated automorphism $\alpha$
of the $C^*$-algebra $C(X)$ of complex-valued functions on $X$ given
by $\alpha (f) = f\circ T$ for all $f\in C(X)$. 
Define the function $f\in C(X)$ by
$$ f((a_k )_{k\in\mathbb{Z}} ) = a_0 $$
for all $(a_k )_{k\in\mathbb{Z}} \in X$. Now given any $n\in\mathbb{N}$,
for every $\gamma\in \{ -1,1 \}^{\{ 0, \dots , n-1 \}}$ there is some
$x\in X$ such that, for each $k=0 , \dots , n-1$, the function $\alpha^k (f)$
takes the value $\gamma (k)$ at $x$. This implies that any map sending
the standard basis of $\ell^n_1$ to $\{ f, \alpha (f) , \dots , 
\alpha^{n-1} (f) \}$ extends linearly to an isometric isomorphism of the real 
linear spans. This simple example illustrates that, in general, a high degree 
of dynamical mixing will produce real linear subspaces which are 
isometrically isomorphic to $\ell^n_1$ in a canonical way with respect to the 
iterates of one or more suitably chosen real-valued functions. 

Now suppose that in our shift
example we have a matrix algebra $M_k$ and completely positive 
contractions $\phi : C(X) \to M_k$ and $\psi : M_k \to C(X)$ such that the 
diagram
\begin{gather*}
\xymatrix{
C(X) \ar[dr]_{\phi} 
\ar[rr]^{\id} & & C(X) \\
& M_k \ar[ur]_{\psi}}
\end{gather*}
approximately commutes to within $\delta$ on $\Omega_n = 
\{ f, \alpha (f) , \dots , \alpha^{n-1} (f) \}$. Then this
diagram also approximately commutes to within $\delta$ on 
the entire unit ball of the
real linear span $X$ of $\Omega_n$ by virtue of the fact that $\Omega_n$ forms
a standard basis for a copy of $\ell^n_1$. It follows that if $\delta < 1$
then $\phi (X)$ is
$(1-\delta )^{-1}$-isomorphic to $\ell^n_1$, i.e., the Banach-Mazur distance
$$ d(\phi (X),\ell^n_1 ) = \inf \{ \| \Gamma \| \| \Gamma^{-1} \| : 
\Gamma : \phi (X)\to\ell^n_1 \text{ is an isomorphism} \} $$ 
is at most $(1-\delta)^{-1}$. Now we know that the Voiculescu-Brown 
entropy agrees with the topological entropy on the 
pure state space, whose positive value of $\log 2$ is 
captured in the dynamical mixing that produces $\ell^n_1$ via the iterates
of $f$. Thus, taking a Banach space viewpoint, we might suspect 
that, in general, for a fixed $K \geq 1$, 
the presence of a subspace $K$-isomorphic to $\ell^n_1$ within 
the real linear space of self-adjoint elements of $M_k$ 
implies an asymptotic exponential dependence of $k$ on $n$.
This is indeed the case, as demonstrated by the following key 
proposition, which
I am grateful to Nicole Tomczak-Jaegermann for having communicated to me.
For simplicity, in the proposition statement and thenceforth (with the 
exception of Proposition~\ref{P-Cantor} below) we will take 
our spaces to be over the complex numbers, as 
this will only 
affect our statements up to a fixed isomorphism factor; for example, in the 
above situation the $(1-\delta )^{-1}$-isomorphism between $\phi (X)$ and 
$\ell^n_1$ extends to a $2(1-\delta )^{-1}$-isomorphism between the complex 
linear span of $\phi (X)$ and the complex scalar version of $\ell^n_1$,
where by a {\it $D$-isomorphism} we mean an isomorphism $\Gamma : Y \to Z$
between Banach spaces which satisfies $\| \Gamma \| \| \Gamma^{-1} \| \leq D$. 
 
\begin{proposition}\label{P-ell1}
Let $X$ be an $n$-dimensional subspace of $M_k$ (with the $C^*$-algebra norm)
which is $D$-isomorphic to $\ell^n_1$. Then 
$$ n \leq a D^2 \log k $$
where $a>0$ is a universal constant.
\end{proposition}

\noindent The proof, which can be found in \cite{EID}, is based on a 
comparison of the (Rademacher) type $2$ constants of the spaces involved. 
The required estimate on the type $2$ constant of $M_k$, in particular, can 
be obtained using bounds on the type $2$ constants of the Schatten 
$p$-classes which follow from Tomczak-Jaegermann's work in \cite{TC}.

With Proposition~\ref{P-ell1} at hand we can now make the following 
general statement yielding positive Voiculescu-Brown entropy as a conclusion. 

\begin{proposition}\label{P-pdpe}
Let $\alpha$ be an automorphism of an exact $C^*$-algebra $A$. Suppose 
there exist a finite subset $\Omega\subset A$,
a $D\geq 1$, and subsets $I_n \subset\{ 0 , \dots , n-1 \} \times\Omega$ 
satisfying $\limsup_{n\to\infty} | I_n | /n > 0$
such that for each $n\in\mathbb{N}$
some (equivalently, any) map from the standard basis of $\ell^{I_n}_1$ to 
$\{ \alpha^k (x) : (k,x)\in I_n \}$ 
linearly extends to a $D$-isomorphism. Then $ht(\alpha ) > 0$.
\end{proposition}

\noindent Proposition~\ref{P-pdpe} is a direct consequence of 
Proposition~\ref{P-ell1}, as the latter shows that if 
$$ \limsup_{n\to\infty} | I_n | /n \geq\mu > 0 $$
then for any $0 < \delta < D^{-1}$ we have
$$ ht(\alpha , \Omega , \delta ) \geq \mu a^{-1} 
D^{-2} (1 - D\delta )^2 > 0 . $$

Note that if the hypotheses of Proposition~\ref{P-pdpe} hold for
an automorphism $\alpha$, then they also hold for any automorphism $\beta$
of a $C^*$-algebra $B$ such that there is a surjective $^*$-homomorphism
$\gamma : B \to A$ satisfying $\gamma\circ\beta = \alpha\circ\gamma$ (as is
easily checked using the contractivity of $\gamma$), so 
that every such $C^*$-dynamical extension $\beta$ of $\alpha$ has
positive entropy. We thus obtain some information about the behaviour of 
Voiculescu-Brown entropy under taking extensions, which in general has 
remained a mystery.

As a concrete illustration of Proposition~\ref{P-pdpe}, consider the following 
operator-theoretic examples from \cite{EID}. We start with a sequence
$\gamma\in\{ -1,0,1 \}^{\mathbb{Z}}$ in which every finite
string of $-1$'s and $1$'s occurs. Setting 
$E_i = \{ k\in\mathbb{Z} : \gamma (k) = i \}$ for each $i=-1,0,1$,
we define the operator $x\in\mathcal{B}(\ell_2 (E_{-1} \cup E_1 ))$ 
by specifying
$$ x\xi_k = \gamma (k) \xi_k $$
on the set $\{ \xi_k : k\in E_{-1} \cup E_1 \}$ of
standard basis elements of $\ell_2 (\mathbb{Z})$.
Next take the direct sum of $x$ with any self-adjoint operator of
norm at most one on $\mathcal{B} (\ell_2 (E_0 ))$, which yields
an operator $y$ on $\mathcal{B} (\ell_2 (\mathbb{Z} )) \supset
\mathcal{B}(\ell_2 (E_{-1} \cup E_1 )) \oplus
\mathcal{B} (\ell_2 (E_0 ))$. Now if $\alpha$ is an automorphism of
an exact $C^*$-subalgebra $A$ of $\mathcal{B} (\ell_2 (\mathbb{Z} ))$
which restricts to the inner automorphism of 
$\mathcal{B} (\ell_2 (\mathbb{Z} ))$ arising from the canonical  
bilateral shift on $\ell_2 (\mathbb{Z} )$ as 
applied to $y$ and its iterates, then $ht(\alpha ) > 0$ by 
Proposition~\ref{P-pdpe}, since it is readily seen from our choice of
sequence $\gamma$ that for every $n\in\mathbb{N}$ any  map from the 
standard basis of $\ell^n_1$ to 
$\{ y, \alpha (y) , \dots , \alpha^{n-1} (y) \}$
linearly extends over $\mathbb{R}$ to an isometric isomorphism, and hence over 
$\mathbb{C}$ to a $2$-isomorphism. Without addressing here the general
question of the existence of exact $C^*$-algebras admitting an automorphism
$\alpha$ as above, we point out in particular that if $a$ is a diagonal
operator then we obtain from the shift on $\ell_2 (\mathbb{Z} )$ a 
positive entropy automorphism of the commutative $C^*$-algebra $A$ 
generated by $a$ and its iterates. Note than we can thus 
obtain positive topological entropy without having a topological 
description of the system at hand.

\begin{problem}\label{PM-conv}
Under what conditions does the converse of Proposition~\ref{P-pdpe} hold?
\end{problem}

For homeomorphisms of the Cantor set we can show that
the converse of Proposition~\ref{P-pdpe} is valid, and even in a stronger
form:

\begin{proposition}\label{P-Cantor}
Let $T:X\to X$ be a homeomorphism of the Cantor set. Then 
$h_{\topol}(T) > 0$ if and only if there is a continuous function $f:X\to
\mathbb{R}$ and subsets $I_n \subset \{ 0 , \dots , n-1 \}$ with
$\limsup_{n\to\infty} | I_n | /n > 0$ such that for each $n\in\mathbb{N}$
the set $\{ f\circ T^k : k\in I_n \}$ forms a standard basis for a copy
of $\ell^{I_n}_1$ inside the real Banach space $C(X,\mathbb{R} )$,
i.e., any map from the standard basis elements of $\ell^{I_n}_1$
to $\{ f\circ T^k : k\in I_n \}$ extends linearly over $\mathbb{R}$ to an 
isometric isomorphism.
\end{proposition}

\begin{proof}
In view of Proposition~\ref{P-pdpe} we need only prove the ``only if''
implication. Suppose then that $h_{\topol}(T) > 0$. We will 
begin by showing the existence
of a $2$-element clopen partition $\mathcal{U}$ of $X$ such that
$h_{\topol}(T,\mathcal{U} ) > 0$ (cf.\ the first part of the proof of 
Proposition~1 in \cite{Bl}). Since the topology of $X$ is generated
by the clopen sets there is a finite clopen cover $\mathcal{V} = \{ V_1 ,
\dots V_n \}$ of $X$ such that $h_{\topol}(T,\mathcal{V} ) > 0$. We may
assume that $\mathcal{V}$ is in fact a partition of $X$ by suitably refining 
it if necessary. For each $k=1, \dots , n$ denote by $\mathcal{V}_k$ the 
clopen partition $\{ V_k , X\setminus V_k \}$. 
Since $\mathcal{V}_1 \vee\dots\vee
\mathcal{V}_n$ is a refinement of $\mathcal{V}$ we then have
$$ 0 < h_{\topol}(T,\mathcal{V} ) \leq h_{\topol} (T, \mathcal{V}_1 \vee\dots
\vee\mathcal{V}_n ) \leq \sum_{k=1}^n h_{\topol}(T,\mathcal{V}_k ) . $$
Thus for some $k=1, \dots , n$ the $2$-element clopen partition
$\mathcal{V}_k$ satisfies $h_{\topol}(T,\mathcal{V}_k ) > 0$, as desired.
Rewrite this clopen partition as $\mathcal{U} = \{ U_1 , U_{-1} \}$.

Next set $f = \chi_{U_1} - \chi_{U_{-1}} \in C(X,\mathbb{R} )$ where
$\chi_{U_1}$ and $\chi_{U_{-1}}$ are the characteristic functions of 
$U_1$ and $U_{-1}$, respectively. For each $n\in\mathbb{N}$ denote by 
$\mathcal{W}_n$ the clopen partition 
$\mathcal{U}\vee T^{-1}\mathcal{U}\vee\cdots\vee T^{-(n-1)} \mathcal{U}$
and by $E_n$ the set of all $\gamma\in
\{ -1,1 \}^{\{ 0, \dots , n-1 \}}$ such that $\bigcap_{k=0}^{n-1}
T^{-k} U_{\gamma (k)} \in\mathcal{W}_n$. By the Sauer-Perles-Shelah lemma 
(or more specifically the consequence thereof formulated as Lemma 2.2 
in \cite{GW}) there are a $c>0$ and an $n_0 \in\mathbb{N}$ 
such that for all $n\geq n_0$ there is a 
subset $I_n \subset\{ 0 , 1 , \dots , n-1 \}$ satisfying $| I_n | \geq cn$
and $E_n | I_n = \{ -1,1 \}^{I_n}$. Now if $n\geq n_0$ then for every
$\gamma\in E_n | I_n$ there is a point $x\in X$ such that
$(f\circ T^k ) (x) = \gamma (k)$ for every $k\in I_n$. As in the
example of the shift on $\{ -1,1 \}^{\mathbb{Z}}$, this implies that
any map sending the standard basis of $\ell^{I_n}_1$ to 
$\{ f\circ T^k : k\in I_n \}$ extends linearly over $\mathbb{R}$ to an
isometric isomorphism,
and thus since $\limsup_{n\to\infty} | I_n | /n \geq c > 0$ we are done.
\end{proof}

In this section we have used only the Banach space structure 
of the $C^*$-algebras in question. For a stark
illustration of how the operator space structure can come into play,
consider the shift 
$u_k \mapsto u_{k+1}$ on the full group $C^*$-algebra 
$C^* (\mathbb{F}_\infty )$ of the free group on a countable set of
generators with associated unitaries $\{ u_k \}_{k\in\mathbb{Z}}$.
For each $n\in\mathbb{N}$ the set $\{ u_1 , u_2 , \dots , u_n \}$ forms a 
standard basis for a copy of $\ell^n_1$, but the $C^*$-algebra
$C^* (\mathbb{F}_\infty )$ is not exact and for $n\geq 2$ any complete
isomorphism between $\spn \{ u_1 , u_2 , \dots , u_n \}$ and a subspace of 
a matrix algebra has completely bounded isomorphism constant at least 
$n(2\sqrt{n-1})^{-1}$ (see \cite{Ex}).
In the next section we will discuss the analogous shift on the reduced
group $C^*$-algebra $C^*_r (\mathbb{F}_\infty )$, which from a geometric 
viewpoint is drastically different.

\section{The free shift}\label{S-free}

Here we apply our geometric perspective to the example of the 
free shift on the reduced group $C^*$-algebra  
$C^*_r (\mathbb{F}_\infty )$ of the
free group on a countable set  $\{ g_k \}_{k\in\mathbb{Z}}$ of generators.
This automorphism, which we will simply refer to as the free shift, 
arises from the shift $k \mapsto k+1$ on the index set $\mathbb{Z}$.
It can be regarded as a quantized version of the shift $k\mapsto k+1$ 
on $\mathbb{Z}$ compactified with a fixed point at infinity, whereby 
orthogonality in $\ell^n_\infty$ (for which the characteristic functions
of the singletons $\{ 1 \} , \dots , \{ n \}\subset\mathbb{Z}$ form 
a standard basis) is replaced by orthogonality 
in $\ell^n_2$ (for which the unitaries associated
to the elements $g_1 , \dots , g_n \in\mathbb{F}_\infty$ form a standard
basis). While this analogy makes little sense from the perspective of 
free probability (in terms of which reduced free products actually exhibit  
parallels with tensor products \cite{NRV}), it is the 
appropriate one for our dynamical context, as we will see. 
In the noncommutative case, however, we must also take the operator
space structure into account. Indeed our 
``quantization'' is not unique and we could also 
consider for example the analogous free shift on the Cuntz 
algebra $\mathcal{O}_\infty$ \cite{BC} (see the second paragraph below).

It was shown in both \cite{TFP} and \cite{BC}, via different arguments,
that the free shift has zero
Voiculescu-Brown entropy. Correspondingly, the compactified shift on
$\mathbb{Z}$ has zero topological entropy, as is readily computed
directly from the definition.

In \cite{NE} St{\o}rmer describes the free shift as the ``most noncommutative''
situation and accordingly asserts that highly 
noncommutative systems tend to have small entropy. The phenomenon 
underlying this statement appears
in fact to have less to do with noncommutativity per se than with the relation 
of orthogonality, which is being considered here in its quantized Hilbert space
sense (cf.\ the discussion in the introduction of \cite{St}) but is equally
well associated to zero entropy (in fact, arithmetic dynamical growth) in
the classical commutative situation. It is hardly a coincidence that 
the matrix algebra $M_n$ accommodates both $\ell^n_\infty$
(down the diagonal) and $\ell^n_2$ (across any row or down any column),
giving us a hint that the dynamical growth is subexponential,
and hence that the entropy is zero, for both 
the compactified shift on $\mathbb{Z}$ and the free shift. This hint
leads directly to a proof of zero Voiculescu-Brown entropy for the 
compactified shift on $\mathbb{Z}$ at the $C^*$-algebra level. On the other
hand, in the noncommutative case we have to be more careful as a 
result of the operator space structures involved, and indeed the 
situation for the free shift 
is subtler and more sophisticated arguments are required \cite{TFP,BC}.
In fact the closed subspace spanned by the unitaries 
in $C^*_r (\mathbb{F}_\infty )$ 
corresponding to the generators of $\mathbb{F}_\infty$ is completely
isomorphic to the closed subspace 
spanned by the elements
$e_{1i} \oplus e_{i1}$  in the direct sum $R\oplus C$ of the row and
column Hilbert spaces in $\mathcal{B} (\ell_2 )$ \cite{HPB} (see also
Section 8.3 of \cite{Pis}). Actually, if we take the operator space 
structure of a matrix algebra into consideration then
the geometric hint from above applies  
precisely and directly in the noncommutative case if we switch to
the nuclear setting of the Cuntz algebra $\mathcal{O}_\infty$ (which can
be viewed as an infinite reduced free product of Toeplitz algebras---see
Chapter 1 of \cite{VDN}) and consider
the automorphism defined by shifting the index on the canonical isometries 
$\{ s_k \}_{k\in\mathbb{Z}}$ \cite{BC}, for then arithmetic dynamical
growth and hence zero Voiculescu-Brown entropy, at least at the local level
of the canonical isometries, is readily seen by combining the fact that 
the closed subspace spanned by $\{ s_k \}_{k\in\mathbb{Z}}$ is canonically
completely isometric to the column Hilbert space in $\mathcal{B} 
(\ell_2 )$ (see Section 1 of \cite{OH})  
with a result of Pop and Smith that permits us to 
use general completely contractive linear maps in the definition of
Voiculescu-Brown entropy \cite{PS} and an appeal to Wittstock's extension
theorem which permits us to extend completely contractive
linear maps into $\mathcal{B} (\mathcal{H})$ for any Hilbert space
$\mathcal{H}$ (see \cite{Pau}).

St{\o}rmer furthermore expresses surprise in \cite{NE} that the free shift has 
zero entropy in view of the fact that it is extremely ergodic, in the sense 
that its extension to the weak operator closure admits no proper globally 
invariant injective von Neumann subalgebra except for the scalars \cite{Po}.
However, extreme ergodicity in topological dynamics (i.e., the non-existence,
at the function level, of proper globally invariant unital
$C^*$-subalgebras besides the scalars, which is usually referred to as
``primeness'' \cite{Prime}) is typically associated
with a lack of recurrence and, in particular, zero entropy (by \cite{Lin} 
all finite-dimensional prime systems have zero entropy). Moreover this
behaviour is manifested in an extreme way in our 
classical example of the compactified shift on $\mathbb{Z}$. Note that 
if $f\in C(\mathbb{Z} \cup \{ \infty \} )$ belongs to a proper unital 
$C^*$-subalgebra of $C(\mathbb{Z} \cup \{ \infty \} )$ which is globally
invariant under the induced $C^*$-algebra automorphism, then with $T$ 
denoting the compactified shift we obtain from the Stone-Weierstrass theorem 
the existence of two distinct points $m,n\in
\mathbb{Z} \cup \{ \infty \}$ such that $f(T^j (m)) = f(T^j (n))$ for
all $j\in\mathbb{Z}$, from which 
we infer by the continuity of $f$ that $f(k) = f(\infty )$ for 
all $k\in\mathbb{Z}$, i.e., $f$ is constant. Thus the only globally invariant
unital $C^*$-subalgebras are the scalars and 
$C(\mathbb{Z} \cup \{ \infty \} )$ itself.

By viewing the free shift as a quantization of the compactified shift
on $\mathbb{Z}$ it therefore seems natural to expect zero Voiculescu-Brown 
entropy in light of the above geometric and topological considerations.
We will not delve here into a rigorous explanation for 
zero entropy, referring again to \cite{TFP,BC} for proofs.

Along the same lines, we might also 
regard a reduced free product of 
$C^*$-algebra automorphisms as a quantized version of a disjoint union 
of homeomorphisms or of a direct sum of $C^*$-algebra automorphisms.
Indeed by \cite{BDS} we have the formula
$$ ht(\alpha * \beta ) = \max (ht(\alpha ), ht(\beta )) $$
for reduced free products with amalgamation over a 
finite-dimensional $C^*$-algebra, paralleling the behaviour of disjoint 
unions with respect to topological entropy or direct sums with respect to
Voiculescu-Brown entropy.

Thus the world of freeness, while inexhaustibly rich from a free 
probability viewpoint (see \cite{VDN,HP}), is associated with a high 
degree of determinism in $C^*$-dynamics.

The discussion of this section vividly 
illustrates the idea, promoted by Weaver in 
\cite{MQ}, that the notion of quantization is at essence about Hilbert space, 
with noncommutativity appearing as a technical consequence.

\section{Induced dynamics on state spaces}\label{S-EID}

We will indicate in this section how Proposition~\ref{P-ell1} can be applied 
in a more systematic way to obtain a noncommutative analogue of Glasner
and Weiss's result (Theorem~\ref{T-GW}) that zero topological entropy 
implies zero entropy on the space of probability measures. Full details 
can be found in \cite{EID}.
 
The proof of Proposition 2.1 in \cite{GW},
which involves a functional-analytic application of the combinatorial
Sauer-Perles-Shelah lemma, can be adapted to establish the following
result. Here $C_1^n$ denotes the space of $n\times n$ matrices
with the trace class norm.

\begin{lemma}\label{L-geom}
Given $\varepsilon > 0$ and $\lambda > 0$ there exist $n_0 \in\mathbb{N}$
and $\mu > 0$ such that, for all $n\geq n_0$, if $\phi : C_1^{r_n} \to
\ell_\infty^n$ is a $^*$-linear map of norm at most $1$ such that the
image of the unit ball of $C_1^{r_n}$ under $\phi$ contains an 
$\varepsilon$-separated set of
self-adjoint elements of cardinality at least $e^{\lambda n}$, 
then $r_n \geq e^{\mu n}$. 
\end{lemma}

Whereas Glasner and Weiss apply information about the possible size of 
subspaces of $\ell^n_\infty$ and $\ell^n_1$ which are
approximately isometric to 
$\ell^k_2$, for the proof of Lemma~\ref{L-geom} we must alternatively appeal 
to Proposition~\ref{P-ell1}, which addresses the matrix situation.   

Now let $A$ be an separable exact $C^*$-algebra and $\alpha$ an automorphism
of $A$, and denote by $T_\alpha$ the induced homeomorphism $\omega\mapsto
\omega\circ\alpha$ of the quasi-state
space $Q(A)$ (i.e., the convex set of positive linear functionals $\phi$ on 
$A$ with $\| \phi \| \leq 1$, equipped with the weak$^*$ topology, under
which it is compact). Suppose that $T_\alpha$ has positive topological
entropy (this actually implies that the entropy is 
infinite---see \cite{EID}---but this is not of consequence for the present
discussion). Then by metrizing $Q(A)$ by taking the supremum of the
absolute values of the differences of two given
elements under evaluation on a compact and total subset $K$ of 
self-adjoint elements of $A$, 
we can apply the separated set definition of topological entropy 
and push everything down to the level of the matrix algebras involved in 
the definition of Voiculescu-Brown entropy to permit an application
of Lemma~\ref{L-geom}, with the required $^*$-linear map from $C_1^{r_n}$
to $\ell^n_\infty$ 
constructed by evaluating a suitable finite subset of $K$ on the relevant 
elements of $Q(A)$ as modeled at the matrix level. Lemma~\ref{L-geom} 
guarantees exponential growth in the rank of the matrix algebras, and so we 
obtain the following result. 

\begin{theorem}\label{T-zero}
Let $A$ be a separable exact $C^*$-algebra and $\alpha$ an
automorphism of $A$. If $\alpha$ has zero Voiculescu-Brown entropy, then
$T_\alpha$ has zero topological entropy.
\end{theorem}

If $A$ is unital then in the statement of the thereom we can also replace 
$T_\alpha$ with its restriction to the state space.

Since topological entropy is non-decreasing under passing to subsystems,
Theorem~\ref{T-zero} shows in particular that if the induced homeomorphim
$T_\alpha$ on the quasi-state space has positive topological entropy, then 
every $C^*$-dynamical extension of $\alpha$ has positive Voiculescu-Brown
entropy (cf.\ the second paragraph after the statement of 
Proposition~\ref{P-pdpe}). More generally:

\begin{corollary}\label{C-pos}
Let $A$ and $B$ be separable exact $C^*$-algebras and 
$\alpha :A\to A$ and
$\beta :B\to B$ automorphisms. Suppose that the homeomorphism of the
quasi-state space of $A$ has positive topological entropy, and 
suppose that there exists a surjective positive contractive linear map
$\gamma : B\to A$ such that $\alpha\circ\gamma = \gamma\circ\beta$, or
an injective positive contractive linear map $\rho : A\to B$ such that
$\beta\circ\rho = \rho\circ\alpha$. Then $\beta$ has positive
Voiculescu-Brown entropy. More generally, the same conclusion
holds whenever $\alpha$ can be obtained from $\beta$ through
a finite chain of intermediary automorphisms intertwined in succession by 
maps of the same form as $\gamma$ or $\rho$.
\end{corollary}

In the other direction, we can draw from Theorem~\ref{T-zero} some
topological-dynamical conclusions such as the following proposition, which
holds in view of the fact that the free shift on
$C^*_r (\mathbb{F}_\infty )$ (see Section~\ref{S-free}) has zero 
Voiculescu-Brown entropy \cite{TFP,BC}.

\begin{proposition}\label{P-fs}
The homeomorphism of the state space of $C^*_r (\mathbb{F}_\infty )$ 
induced by the free shift has zero topological entropy.
\end{proposition}

Finally, we ask:

\begin{question}
Does the converse of Theorem~\ref{T-zero} hold?
\end{question}

\noindent Possible candidates for counterexamples to the converse of 
Theorem~\ref{T-zero} can be found among the collection of automorphisms 
of the rotation $C^*$-algebras $A_\theta$ associated to a given matrix
$S\in SL(2,\mathbb{Z})$ with eigenvalues off the unit circle. For every
$\theta$ the automorphism of $A_\theta$ defined via $S$ has positive 
Voiculescu-Brown entropy \cite{NTA}, but we have not been able to
determine the topological entropy on the state space for irrational
$\theta$.
See \cite{EID} for details.
\bigskip

\noindent{\it Note added:} An answer to Problem~\ref{PM-conv} has now been
given in \cite{PTE}.

\end{document}